\documentclass{amsart}

\usepackage{amssymb}
\usepackage{graphicx}
\usepackage[T1]{fontenc}
\usepackage{hyperref}

\newtheorem{theorem}{Theorem}[section]
\newtheorem{lemma}[theorem]{Lemma}

\theoremstyle{definition}

\theoremstyle{remark}
\newtheorem{remark}[theorem]{Remark}

\numberwithin{equation}{section}

\begin{document}

\title{Infinite family of elliptic curves of rank at least 4}
\author{Bartosz Naskr\k{e}cki}
\address{Graduate School, Department of
Mathematics and Computer Science, Adam Mickiewicz University, Pozna\'{n}, Poland}
\email{nasqret@gmail.com}


\date{\today}


\begin{abstract}
We investigate $\mathbb{Q}$-ranks of the elliptic curve $E_t$: $y^2+txy=x^3+tx^2-x+1$ where $t$ is a rational parameter. We prove that for infinitely many values of $t$ the rank of $E_t(\mathbb{Q})$ is at least 4.
\end{abstract}

\maketitle

\section{Introduction}
It is of fundamental interest to find families of elliptic curves parametrized by a rational parameter with ranks higher than a prescribed constant cf. \cite{3}. In this paper we investigate the family of curves:
\begin{equation}\label{curve0}
E_t:y^2+txy=x^3+tx^2-x+1
\end{equation}
with parameter $t\in\mathbb{Q}$. We prove the following
\begin{theorem}\label{main_theorem}
For infinitely many $u\in\mathbb{Q}$ the elliptic curve over $\mathbb{Q}$ given by the affine equation $E_{(u^2-u-3)}:y^2+(u^2-u-3)xy=x^3+(u^2-u-3)x^2-x+1$ has the Mordell-Weil group of rank at least 4. More precisely, the group $E_{(u^2-u-3)}(\mathbb{Q})$ contains the subgroup spanned by the following linearly independent points:

$$(0,1),(1,1),(u,u+1),\left(\frac{1}{9},\frac{1}{54} (9 + 3 u - 3 u^2 + v)\right)$$
where the point $(u,v)$ lies on the elliptic curve given by the equation:
$$2569 + 18 u - 9 u^2 - 18 u^3 + 9 u^4=v^2.$$
The latter curve has the Weierstrass model:
\begin{equation}\label{explicit_eq}
y_0^2 = x_0^3 - 92835 x_0 + 1389150
\end{equation}
which defines elliptic curve over $\mathbb{Q}$ with the Mordell-Weil group of rank 2, spanned by the points: 
$$(-309, 756),(-45,2340),(15,0)$$
\end{theorem}
We have checked that the curves $E_{(u^2-u-3)}$ have different j-invariant for all but finitely many $u\in\mathbb{Q}$ as in the statement of the theorem.
\newline
\newline
\noindent
In \cite{1} Brown and Myers constructed an infinite family of elliptic curves over $\mathbb{Q}$ with quadratic growth of parameter and the rank of the Mordell-Weil group at least 3. They proved the linear independence of classes in $E(\mathbb{Q})/2E(\mathbb{Q})$ for curves $E$ with trivial torsion. We investigate the question posed at the end of the paper \cite{1} to find families of elliptic curves with the rank of the Mordell-Weil group higher than 3.

With the family defined in Theorem \ref{main_theorem} we obtained rank 4 at cost of using Mordell's Conjecture,i.e. theorem of Faltings.
In \cite{2} Rubin and Silverberg also obtained an infinite family of twists of elliptic curves over $\mathbb{Q}$ of rank 4. The method of \cite{2} relies on Silverman's specialization theorem cf. \cite{5}. The investigated curves are in Legendre's form. The families are parametrized by projective line or by elliptic curve of rank 1. The twists are parametrized by another elliptic curve of rank 1.

In Theorem \ref{main_theorem} the parametrization of the family is given by the elliptic curve of rank 2 and we use a single parameter instead of two. Our methods are more elementary and in some places more explicit, since we don't use specialization theorems and the coefficients of the family are at most quadratic. 
\begin{remark}
We can extend the result stated in Theorem \ref{main_theorem} following \cite{4}.  Let $\Lambda(E/\mathbb{Q},s)$ be the complete L-series of the elliptic curve over $\mathbb{Q}$. Denote by $W(E)\in\{\pm 1\}$ the root number in the functional equation:
\begin{equation}
\Lambda(E,2-s)=W(E)\Lambda(E,s)
\end{equation}
The Parity Conjecture predicts that:
\begin{equation}
(-1)^{\textrm{rank}E(\mathbb{Q})}=W(E)
\end{equation}
In our particular case we can compute the root number $W(E_u)$ for the specific curves $E_u$  and determine the parity of the rank of group $E_u(\mathbb{Q})$. Since the computations can be done explicitly using primes of bad reduction of $E_u$ (for example using SAGE) we state the numerical results in Section \ref{num_results}. Assuming Parity Conjecture we construct several elliptic curves over $\mathbb{Q}$ that have rank at least 5 (Table \ref{tab:5}).
\end{remark}
\section{Description of the algorithm}
There are two obvious points lying on the curve (\ref{curve0}), namely:
$$(0,1),(1,1)\in E_u(\mathbb{Q}(u))$$
We produce with several points with coordinates in the ring $\mathbb{Z}[u]$:
\begin{eqnarray}
-(0,1)& =& (0,-1) \nonumber\\
-(1,1)& =& (1,-u-1) \nonumber\\
(0,1)+2(1,1)& =& (-u + 1,-1)\nonumber\\
(0,1)+(1,1)& = & (-u - 1, u^2 + u - 1)\nonumber\\
(0,1)-(1,1)& = & (u + 3, 2u + 5)\nonumber\\
-(0,1)+(1,1)& = & (u + 3, -u^2 - 5u - 5 )\nonumber\\
-(0,1)+2(1,1)& = & (u + 5, 2u + 11)\nonumber\\
2(1,1)& = & (-1, u + 1)\nonumber
\end{eqnarray}
In order to find more points on the curve (\ref{curve0}) we specialize parameter $u$ to a polynomial function of another parameter $v\in\mathbb{Q}$:
$$u(v)=a_n v^n +\ldots+a_1 v+a_0$$
where $a_i\in\mathbb{Q}$. To get a rational point on the curve (\ref{curve0}) with x-coordinate $av+b$ and $a,b\in\mathbb{Q}$ it is necessary and sufficient:
\begin{equation}\label{discriminant1}
\Delta(v)=4(-1+b+a v)^2 (1+b+a v)+(2+u(v))^2 (b+a v)^2
\end{equation}
is a perfect square.
\begin{lemma}\label{coeff}
Let $P=(x,y)$ be a rational point on the curve $E_u(v)$ over $\mathbb{Q}(v)$ where $u(v)=a_n v^n +\ldots+a_1 v+a_0\in\mathbb{Q}[v]$ of positive degree.
Let $x=a v+b\in\mathbb{Q}[v]$. Since $a\neq 0$ it follows that:
\begin{enumerate}
    \item[(i)] $deg(u)=1\Longrightarrow $ $P = k (0,1)+l(1,1)$ for some $k,l\in\mathbb{Z}$
    \item[(ii)] $deg(u)=2\Longrightarrow $ $P=(x,x+1)$ and $u=x^2-x-3$ or $P= (x,x-1)$ and $u=-x^2+x+1$
    \item[(iii)] For any $u(v)$ with $deg(u)>2$ there is no rational point with x-coordinate equal to $a v+b\in\mathbb{Q}[v]$ with $a\neq 0$.
\end{enumerate}
\end{lemma}
\begin{proof} 
Assume $u=u(v)=a_1 v+a_0$ and $a_1,a\neq 0$. Put $\Delta(v)=(q_2 v^2+q_1 v+q_0)^2$ and comparing coefficients of (\ref{discriminant1}) and $P(v)$ in descending order we find:
$$q_2=\epsilon a a_1$$
$$q_1=\frac{2 a^2+2 a a_1+a a_0 a_1+b a_1^2}{\epsilon a_1}$$
$$q_0=\frac{-2 a^3-4 a^2 a_1-2 a^2 a_0 a_1-2 a a_1^2+4 a b a_1^2+2 b a_1^3+b a_0 a_1^3}{\epsilon  a_1^3}$$
for $\epsilon=\pm 1$. Equating the last two coefficients gives two equations in variables $a_1,a_0,a,b$:
\begin{align*}
R_1(a,b,a_0,a_1)=-a^6-4 a^5 a_1-2 a^5 a_0 a_1-6 a^4 a_1^2+4 a^4 b a_1^2-4 a^4 a_0 a_1^2\\-a^4 a_0^2 a_1^2-4 a^3 a_1^3+10 a^3 b a_1^3-2 a^3 a_0 a_1^3+5 a^3 b a_0 a_1^3-a^2 a_1^4\\
+8 a^2 b a_1^4-4 a^2 b^2 a_1^4+4 a^2 b a_0 a_1^4+a^2 b a_0^2 a_1^4\\+2 a b a_1^5-4 a b^2 a_1^5+a b a_0 a_1^5-2 a b^2 a_0 a_1^5+a_1^6-b a_1^6-b^2 a_1^6+b^3 a_1^6=0
\end{align*}
\begin{align*}
R_2(a,b,a_0,a_1)=2 a^4+6 a^3 a_1+3 a^3 a_0 a_1+6 a^2 a_1^2-3 a^2 b a_1^2+4 a^2 a_0 a_1^2\\+a^2 a_0^2 a_1^2+2 a a_1^3-4 a b a_1^3+a a_0 a_1^3-2 a b a_0 a_1^3-a_1^4-b a_1^4+b^2 a_1^4=0
\end{align*}
The ideal $I=I(R_1,R_2)$ of the equations can be rearranged to the form of the Gr\"{o}bner basis $I=I(a^9-2 a^7 a_1^2+a^5 a_1^4,R'_1,\ldots,R'_{18})$ with $R'_i=R'_i(a,b,a_0,a_1)$. It implies $a=\pm a_1$. For $a=a_1$ the equations reduce to:
$$a_0=-3+b$$
or
$$a_0=-5+b.$$
For $t=av+b-3$ or $t=av+b-5$ we get points $(t+3,5 + 2 t),(t+3,-5 - 5 t - t^2)$ and $(t+5,11 + 2t),(t+5,-11 - 7 t - t^2)$ respectively. If $u=t$ all these points are linear combinations of $(0,1)$ and $(1,1)$ as on the list above.
For $a=-a_1$ the equations reduce to:
$$a_0=-1-b$$
or
$$a_0=1-b.$$
For $t=-av - b - 1$ or $t=-av - b + 1$ we get points $(-t - 1,1),(-t - 1,-1 + t + t^2)$ and $(-t + 1,-1),(-t+1,1 - t + t^2)$ respectively. Again for $u=t$ all points are the linear combinations of $(0,1)$ and $(1,1)$.
\medskip

\noindent
Assume $u=u(v)=a_2 v^2+a_1 v+a_0$ and $a_1,a\neq 0$. Put $P(v)=(q_3 v^3+q_2 v^2+q_1 v+q_0)^2$. Comparing coefficients of (\ref{discriminant1}) and $P(v)$ implies:
$$q_3=\epsilon  a a_2$$
$$q_2= \frac{a a_1+b a_2}{\epsilon }$$
$$q_1=\frac{2 a+a a_0+b a_1}{\epsilon }$$
$$q_0=\frac{2 a^2+2 b a_2+b a_0 a_2}{\epsilon  a_2}$$
with $\epsilon=\pm 1$.
This gives:
$$a_1=\frac{(-1+2 b) a_2}{a}$$
$$a^3 \left(2+a_0\right)+a \left(1+b-b^2\right) a_2=0$$
$$a^2=\lambda a_2$$
with $\lambda=\pm 1$. It implies $u(v) = -2-\lambda -t \lambda +t^2 \lambda $, where $t=a v+b$. In this way we get two distinct families:
\begin{table}[h]
\caption{Families}
\label{tab:1}
\centering
\begin{tabular}{|c|c|c|}
\hline
 & \textbf{Point} & \textbf{Parameter u} \\
\hline
Family A & $(t,t+1)$ & $t^2-t-3$\\
\hline
Family B & $(t,t-1)$ & $-t^2+t+1$\\
\hline
\end{tabular}
\end{table}
\medskip

\noindent
Consider $u(v)$ as a polynomial in $v$ of degree $n>2$ then $deg(\Delta)=2n+2$ so we look for the polynomial $P(v)$ of degree $n+1$.
We put:
$$a^{*}_{i}=a_i$$
$$q^{*}_{j}=q_j$$
for $0< i\leq n$ and $0\leq j\leq n+1$,$a^{*}_{0}=a_0+2$. For other $i,j$ put:
$$a^{*}_{i}=0$$
$$q^{*}_{j}=0$$
We prove by induction the following formula:
$$q^{*}_{j}=\epsilon(a a^{*}_{j-1}+b b^{*}_j)$$
using the identities:
$$c_j(\Delta)=a^2\sum_{j=\alpha+\beta}a^{*}_{\alpha}a^{*}_{\beta}+2ab\sum_{j+1=\alpha+\beta}a^{*}_{\alpha}a^{*}_{\beta}+b^2\sum_{j+2=\alpha+\beta}a^{*}_{\alpha}a^{*}_{\beta}$$
$$c_j(P^2)=\sum_{j=\alpha+\beta}q^{*}_{\alpha}q^{*}_{\beta}$$
for $j=n+1,...,2n+2$, where $c_j(a_0+a_1 x+\ldots+a_n x^n) = a_j$.It follows from  $\Delta=P^2$:
$$c_j(\Delta)=c_j(P^2).$$
We substitute the coefficients $q_0=\epsilon(2b+b a_0)$ and $q_1=\epsilon (a(a_0+2)+b a_1)$ into the identities above with $j=0,1,2$ and we get $b^2=1$ and finally $a=0$, a contradiction. This completes the proof of the lemma.
\end{proof}
By Lemma \ref{coeff}, we can specialize to one of the quadratic parameters, since the families given by Table \ref{tab:1} have similar properties.
By abuse of notation, we use the same letter $u^2-u-3$ for the parameter u. From now on the curve $E_u$ will be defined by the equation: $$y^2+(u^2-u-3)xy=x^3+(u^2-u-3)x^2-x+1.$$
For simplicity of notation, we write $f(x,y,u)=0$ instead of the equation above.
The point $(u,u+1)$ lies on these curves and gives several new integral points over $\mathbb{Q}[u]$:
\begin{eqnarray}
-(u,u+1)& =& (u, -u^3 + u^2 + 2u - 1)  \nonumber\\
(0,1)+(u,u+1)& =& (-u + 1, u^3 - 2u^2 - u + 1) \nonumber\\
(1,1)-(u,u+1)& =& (u^3 - 2u, u^4 + u^3 - 3u^2 - 2u + 1)\nonumber\\
2(1,1)+(u,u+1)& = & (-u^3 + 4u^2 - 6u + 4, u^5 - 6u^4 + 14u^3 - 17u^2 + 10u - 1) \nonumber
\end{eqnarray}
In order to find the fourth linearly independent rational point on the curve $E_u$ we consider the following general algorithm:
\begin{enumerate}
    \item Choose two rational functions $a(x),b(x)\in\mathbb{Q}(x)$.
    \item Set the simultaneous equations of the form:
\begin{eqnarray}
f(a(u),y_a(u),u)&=&0\nonumber\\
f(b(u),y_b(u),u)&=&0\nonumber
\end{eqnarray}
    \item Find such $a(x),b(x)$ that $y_a(x),y_b(x)\in\mathbb{Q}(x)$.
    \item Sufficient and necessary condition for $y_a,y_b$ to be rational is that the discriminant of the quadratic equation $f(a(x),y_a,x)=0$ in $y_a$  is a perfect square. The same condition holds for the equation in $y_b$.
    \item Find all rational points (the triples $(u,s,t)\in\mathbb{Q}^3$ on the affine curve):
\begin{eqnarray}
\label{curve2}\Delta_{f(a(x),y_a,x)=0}(u)&=&s^2\\
\label{curve3}\Delta_{f(b(x),y_b,x)=0}(u)&=&t^2\nonumber
\end{eqnarray}
\noindent
where $\Delta_{f(a(x),y_a,x)=0}(x),\Delta_{f(b(x),y_b,x)=0}(x)\in\mathbb{Q}[x]$
\end{enumerate}
We pick now $a(x)=x$ and $b(x)=c$. The first equation (\ref{curve2}) reduces to:
$$(2 - u - u^2 + u^3)^2=s^2$$
while the second gives:
$$4 - 4 c - 3 c^2 + 4 c^3 + 2 c^2 u - c^2 u^2 - 2 c^2 u^3 + c^2 u^4=t^2.$$
We choose such a $c\in\mathbb{Q}$ that it defines the elliptic curve in quartic form with infinitely many points $(u,t)$.
Direct search with $u\in\mathbb{N}$ reveals that for $u=7$ we have on the curve $E_7$ four linearly independent points:
$$(0,1),(1,1),(7,8),(\frac{1}{9},\frac{8}{27})$$
though we put $c=\frac{1}{9}$(as in the statement of Theorem \ref{main_theorem}).

\section{Proofs}
In order to prove Theorem \ref{main_theorem} we will need the following elementary lemma:
\begin{lemma}\label{BM}
Let M be a left $\mathbb{Z}-module$. Given any $b\in M$ suppose that $a_1,\ldots,a_k\in M$ are linearly independent over $\mathbb{Z}$ and the nonzero cosets $[a_1],\ldots,[a_k] \in M/2M$ are linearly independent over $\mathbb{F}_2$. 
If $[b]\notin \langle [a_1],\ldots,[a_k]\rangle$ and $M[2]={0}$, two-torsion of M is trivial, then $b,a_1,\ldots,a_k$ are independent over $\mathbb{Z}$ in M.
\end{lemma}
\begin{proof}
Suppose, contrary to our claim, that there exists $\alpha_1,\ldots,\alpha_k,\beta\in\mathbb{Z}$,not all zero, such that:
$$\beta b+\alpha_1 a_1+\ldots +\alpha_k a_k=0.$$
We can assume that $\beta$ is the least positive integer for which the last equation holds.
If $\beta$ is odd:
$$[\beta b]=[b]$$
and
$$[b]=[\alpha_1 a_1+\ldots +\alpha_k a_k]$$
which is a contradiction.
If $\beta$ is even:
$$[0]=[\alpha_1 a_1+\ldots +\alpha_k a_k]$$
The linear independence of cosets $[a_i]$ over $\mathbb{F}_2$ implies that all $\alpha_i$ are even. So it is possible to write:
$$\beta ' b=\alpha_1 ' a_1+\ldots+\alpha_k ' a_k$$
where $2\beta '=\beta$ and $2\alpha_i '=\alpha_i$. Again the contradiction with the minimality of $\beta$.
\end{proof}
Let $M=E_u(\mathbb{Q})$. First we compute the rational 2-torsion for the curve $E_u$. For a point $P=(x,y)$ on the curve $y^2+(u^2-u-3)xy=x^3+(u^2-u-3)x^2-x+1$ the negative is equal:
$$-P=(x,(-u^2 + u + 3)x - y)$$
We find the following condition for the point to be a 2-torsion point:
$$1 - x + \frac{1}{4}(4 (-3 - u + u^2) + (-3 - u + u^2)^2) x^2 + x^3=0$$\label{curve1}
Computing with Maple we find that the equation defines a hyperelliptic curve of genus 2. With substitutions:
\begin{align}
x_0& = \frac{-1-ux+x u^2}{x-1}\\
y_0& = -4x+8ux
\end{align}
the curve transforms birationally into the form:
\begin{equation}\label{eq1002}
y_0^2=-55 + 192 x_0 - 114 x_0^2 + 68 x_0^3 - 15 x_0^4 + 4 x_0^5.
\end{equation}
\begin{lemma}
For all but finitely many $u\in\mathbb{Q}$ the rational 2-torsion of elliptic curve $E_u$ over $\mathbb{Q}$ is trivial.
\end{lemma}
\begin{proof}
This is an immediate corollary of the above calculation and Mordell's conjecture proven by Faltings, since the curve (\ref{eq1002}) has nonsingular model and the genus is 2.
\end{proof}
\begin{proof}[Proof of Theorem \ref{main_theorem}]
In order to prove Theorem \ref{main_theorem} we check if the apropriate points and their linear combinations belong to $2E_u(\mathbb{Q})$. Given a $\mathbb{Q}$-rational point $P=(x,y)$ on the curve $E_u$ over $\mathbb{Q}$ we have the following formula for the x-coordinate of the point $2P$:
$$x(2P)=\frac{4 - 2 u + u^2 + 2 u^3 - u^4 - 8 x + 2 x^2 + x^4}{4 - 
 4 x + (-3 + 2 u - u^2 - 2 u^3 + u^4) x^2 + 4 x^3}.$$
To ease notation define:
$$P_{\varepsilon_1,\varepsilon_2,\varepsilon_3}=\varepsilon_1
(0,1)+\varepsilon_2 (1,1)+\varepsilon_3 (u,u+1)$$
If for $u\in\mathbb{Q}$ there exists a rational point $(\frac{1}{9},y)$ on the curve $E_u$ and $y$ is one of two values:
$$y=\frac{1}{54} (9 + 3 u - 3 u^2 \pm \sqrt{2569 + 18 u - 9 u^2 - 18 u^3 + 9 u^4})$$
then we put:
$$Q_{\varepsilon_1,\varepsilon_2,\varepsilon_3,\varepsilon_4}=\varepsilon_1
(0,1)+\varepsilon_2 (1,1)+\varepsilon_3 (u,u+1)+\varepsilon_4 (1/9,y)$$
where $\varepsilon_i\in\{-1,0,1\}$.
The proof separates naturally into two distinct parts. In the first part we establish the criteria for which the equations:
$$P_{\varepsilon_1,\varepsilon_2,\varepsilon_3}=x(2P)$$
and
$$Q_{\varepsilon_1,\varepsilon_2,\varepsilon_3,\varepsilon_4}=x(2P)$$
have solutions in pairs of rational numbers $(u,x)$ (recall that $P=(x,y)$ lies on $E_u$).
In the second part of the proof we gather information to find the infinite subset of $\mathbb{Q}$ of parameters $u$ for which the rank is at least 4.
To use Lemma \ref{BM} we must consider the tuples:
\begin{equation}
(\varepsilon_1,\varepsilon_2,\varepsilon_3)\in\{(1,0,0),(0,1,0),(0,0,1),(1,1,0),
(0,1,-1),(1,0,1),(1,1,1)\}
\end{equation}
Assume that $(\frac{1}{9},y)$ is $\mathbb{Q}$-rational. We consider:
\begin{multline}
(\varepsilon_1,\varepsilon_2,\varepsilon_3,\varepsilon_4)\in\{(0,0,0,1),(1,0,0,1),(0,1,0,1),\\
(0,0,1,1),(1,1,0,1),(0,1,-1,1),(1,0,1,1),(1,1,1,1)\}
\end{multline}
The tuples with negative entries were chosen to lower the genera of corresponding curves. Since we work mod $2 E(\mathbb{Q})$ the tuples can be chosen quite arbitrarily. We compute genera of curves using \textit{genus} command from \textit{algcurves} package in Maple 12. 
We consider the following separate cases:
\medskip

\noindent
{\it Case I: \rm} 
$(\varepsilon_1,\varepsilon_2,\varepsilon_3)=(1,0,0)$
\newline
The tuple implies the following equation:
\begin{equation}
\frac{4 - 2 u + u^2 + 2 u^3 - u^4 - 8 x + 2 x^2 + x^4}{
 4 - 4 x + (-3 + 2 u - u^2 - 2 u^3 + u^4) x^2 + 4 x^3} = 0
\end{equation}
Since $(0,1)\neq\mathcal{O}$, the point at infinity, then the denominator is nonvanishing and:
$$4 - 2 u + u^2 + 2 u^3 - u^4 - 8 x + 2 x^2 + x^4 = 0$$
It defines an elliptic curve of rank 1. By the following formulas:
\begin{align*}
x_0& =\frac{1}{3 \left(u^2-u-1\right)}(12 u^4-12 u^3 x-36 u^3+12 u^2 x^2+30 u^2 x \\
 & +35 u^2-12 u x^3-24 u x^2-42 u x+61 u+18 x^3+6 x^2+36 x-113)\\
y_0& = -\frac{1}{u^2-u-1}(2 (8 u^5-8 u^4 x-32 u^4+8 u^3 x^2+28 u^3 x\\
 & +44 u^3-8 u^2 x^3-24 u^2 x^2-39 u^2 x+9 u^2+20 u x^3\\
 & +20 u x^2+43 u x-101 u-19 x^3-11 x^2-54 x+122))
\end{align*}
the equation transforms to the short Weierstrass form:
$$y_0^2=x_0^3+\frac{359}{3}x_0+\frac{3130}{27}.$$
The Mordell-Weil group of this elliptic curve is generated by the point $(\frac{53}{3},88)$. Hence in the original form the generator is equal to $(u,x)=(\frac{1}{2},\frac{1}{2})$.
Similar calculations and computation with Maple of the genus are summarized in the table (\ref{tab:2}).
\begin{table}[h]
\caption{Other curves}
\label{tab:2}
\centering
\begin{tabular}{|c|c|}
  \hline 
  \textbf{$(\varepsilon_1,\varepsilon_2,\varepsilon_3)$} & \textbf{genus}\\
  \hline
  $(1,0,0)$ & 1 \\ 
  \hline
  $(0,1,0)$ & 3 \\
  \hline
  $(0,0,1)$ & 2 \\
  \hline
  $(1,1,0)$ & 3 \\
  \hline
  $(0,1,-1)$& 4 \\
  \hline
  $(1,0,1)$ & 2 \\
  \hline
  $(1,1,1)$ & 4 \\
  \hline 
\end{tabular}
\end{table}
Assume now that we are given a point $(\frac{1}{9},y)$ rational, lying on the curve $E_u$ for a suitable $u\in\mathbb{Q}$.
\medskip

\noindent
{\it Case II: \rm} 
$(\varepsilon_1,\varepsilon_2,\varepsilon_3,\varepsilon_4)=(0,0,0,1)$
\begin{equation}
\frac{-u^4+2 u^3+u^2-2 u+x^4+2 x^2-8 x+4}{\left(u^4-2 u^3-u^2+2 u-3\right) x^2+4 x^3-4 x+4}=\frac{1}{9}
\end{equation}
The equation above defines affine curve of genus 5.
\medskip

\noindent
{\it Case III: \rm} 
$(\varepsilon_1,\varepsilon_2,\varepsilon_3,\varepsilon_4)=(1,0,0,1)$
\begin{equation}\label{eq1001}
\frac{-u^4+2 u^3+u^2-2 u+x^4+2 x^2-8 x+4}{\left(u^4-2 u^3-u^2+2 u-3\right) x^2+4 x^3-4 x+4}=-9 u^2+9 u-162 y+180
\end{equation}
From the equation of the curve $E_u$ we can find the following formula:
$$y=\frac{1}{54} \left(-3 u^2+3 u+v+9\right)$$
$$v=\pm \sqrt{2569+18 u-9 u^2-18 u^3+9 u^4}$$
Using these relations we can assume that if a point (x,u) lies on the curve defined as in (\ref{eq1001}) then it also lies on the curve:

\begin{multline}
9(9 u^4-18 u^3-9 u^2+18 u+2569) (u^4 x^2-2 u^3 x^2-u^2 x^2\\
+2 u x^2+4 x^3-3 x^2-4 x+4)^2-(153 u^4 x^2\\
+u^4-306 u^3 x^2-2 u^3-153 u^2 x^2-u^2\\
+306 u x^2+2 u-x^4+612 x^3-461 x^2-604 x+608)^2=0
\end{multline}
This curve has genus 9.
The rest is computed in a similar way (cf. Table \ref{tab:3}).
\begin{table}[h]
\caption{Other curves}
\label{tab:3}
\centering
\begin{tabular}{|c|c|}
  \hline 
  \textbf{$(\varepsilon_1,\varepsilon_2,\varepsilon_3,\varepsilon_4)$} & \textbf{genus}\\
  \hline
  $(0,0,0,1)$ & 5 \\ 
  \hline
  $(1,0,0,1)$ & 9 \\
  \hline
  $(0,1,0,1)$ & 13 \\
  \hline
  $(0,0,1,1)$ & 19 \\
  \hline
  $(1,1,0,1)$& 13 \\
  \hline
  $(0,1,-1,1)$& 15 \\
  \hline
  $(1,0,1,1)$ & 11 \\
  \hline 
  $(1,1,1,1)$ & 28 \\
  \hline
\end{tabular} 
\end{table}
In the last step of the proof we show for which $u\in\mathbb{Q}$ the point $(\frac{1}{9},y)$ is $\mathbb{Q}$-rational.
It is straightforward to compute the root $y$:
$$-\frac{649}{729}+\frac{1}{81} \left(3+u-u^2\right)+\frac{1}{9} \left(-3-u+u^2\right) y+y^2=0$$
$$y=\frac{1}{54} \left(9+3 u-3 u^2\pm\sqrt{2569+18 u-9 u^2-18 u^3+9 u^4}\right)$$
Hence $y\in\mathbb{Q}$ if and only if $2569+18 u-9 u^2-18 u^3+9 u^4$ is a full square. This condition defines the elliptic curve in a quartic form. It is birational to elliptic curve in the Weierstrass form:
$$y_0^2 = x_0^3 - 92835 x_0 + 1389150.$$
The Mordell-Weil group of the curve has rank 2. The torsion subgroup is isomorphic to $\mathbb{Z}/2$. Generators of the free part are $(x_0,y_0)=(-309, -756), (390, -4950 )$ and the generator of the torsion subgroup is $(15,0)$.
In the quartic form they correspond to:
$$\left(\frac{1}{9},\frac{1369}{27}\right),\left(\frac{27}{10},\frac{-5173}{100}\right),\left(-6,-133\right)$$ respectively.

Consider the curve $y^2+uxy=x^3+u x^2-x+1$. The only critical value which gives an infinite subset of parameters $u\in\mathbb{Q}$ was obtained from doubling the point $(0,1)$. We consider this case if the point $(\frac{1}{9},y)$ on the curve $y^2+uxy=x^3+u x^2-x+1$ is $\mathbb{Q}$-rational. Solving the quadratic equation gives:
$$y=\frac{1}{54} \left(-3 u\pm\sqrt{2596+36 u+9 u^2}\right).$$
The point $(0,1)$ is a double and $(\frac{1}{9},y)$ is a rational point when we have the rational point on the curve:
\begin{align}
s^2&= 2596 + 36 u + 9 u^2\\
0& =1 - 4 u - u^2 - 8 x + 2 x^2 + x^4.
\end{align}
Parametrization of the first equation gives:
\begin{align}
u&=\frac{2596-t^2}{6 t-36}\\
s&=\frac{t^2-12 t+2596}{2 t-12}
\end{align}
with a new parameter $t\in\mathbb{Q}\backslash \{6\}$. Substituting into the second equation gives the curve:
\begin{multline}
(-6364096 - 62736 t + 5084 t^2 + 24 t^3 - t^4 - 10368 x \\
+ 3456 t x - 288 t^2 x + 2592 x^2 - 864 t x^2 + 72 t^2 x^2 \\
+ 1296 x^4 - 432 t x^4 + 36 t^2 x^4)
\end{multline}
This curve has genus 3, so it has finitely many rational points. Specializing to a parameter $u^2-u-3$ we obtain that there are only finitely many $u\in\mathbb{Q}$ for which $(0,1)$ is a double while $(\frac{1}{9},y)$ is $\mathbb{Q}$-rational.
From the infinite family of parameters $u$ chosen from the first coordinate of the elliptic curve over $\mathbb{Q}$:
$$2569+18 u-9 u^2-18 u^3+9 u^4=v^2$$
we have excluded only finitely many $u$ due to the doubling restrictions given in Table \ref{tab:2} and Table \ref{tab:3}.
It remains to show that the j-invariant of the curve $F_u: y^2+uxy=x^3+ux^2-x+1$ repeats itself for finitely many $u\in\mathbb{Q}$. We compute:
\begin{equation}
j(F_u)=-\frac{(48 + u^2 (4 + u)^2)^3}{(2 + u)^2 (92 + (-1 + u) u (4 + u) (5 + u))}
\end{equation}
Hence the equation:
$$j(F_u)=j(F_v)$$
defines an affine curve with coordinates $(u,v)$ which has genus 11 according to computations in Maple. This implies that specializing the parameter $u$ to $u^2-u-3$ gives a curve with finitely many points.
\end{proof}

\section{Numerical results}\label{num_results}
\subsection{General statistics}
\begin{figure}[h]
\caption{Curves of rank 3}
\centering
\includegraphics{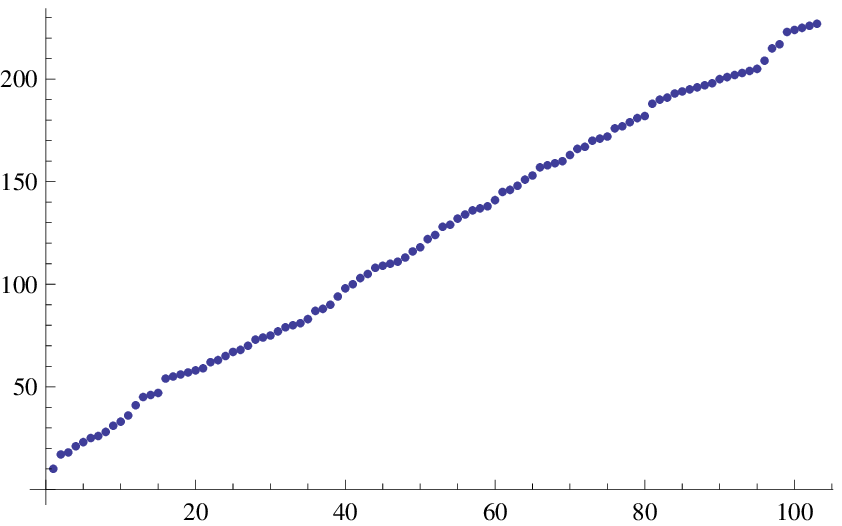}
\label{fig:2}
\end{figure}

We show the statistics of ranks for the family of curves:
$$y^2+uxy=x^3+ux^2-x+1$$
with parameter $u\in\mathbb{N}$. All computations were performed with SAGE 3.4 with mwrank procedure. For some points there is no proof for the upper bound of the rank (in that case the value is excluded).
\begin{figure}[h]
\caption{Rank statistics}
\centering
\includegraphics{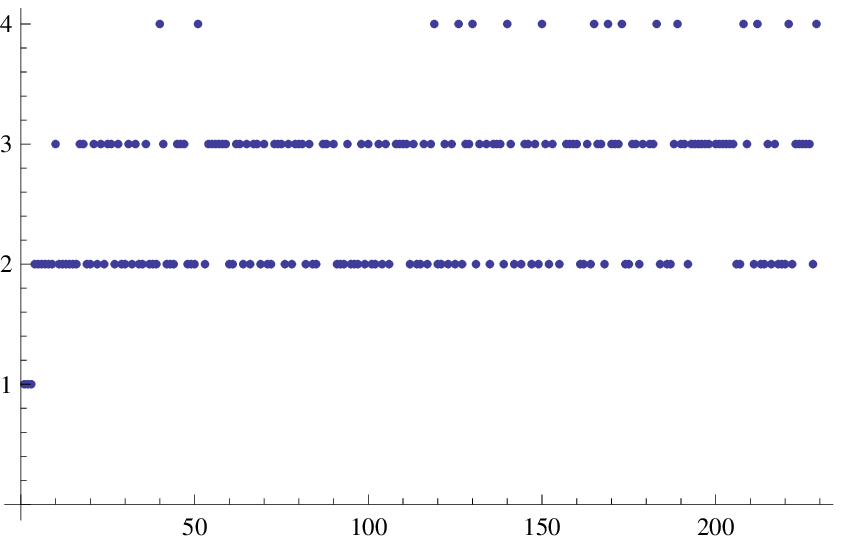}
\label{fig:1}
\end{figure}
We can compute the percentage of curves of each type (Table \ref{tab:4}).
\begin{table}[h]
\caption{Curves of certain ranks}
\centering
\begin{tabular}{|c|c|}
\hline
 \textbf{rank} & \text{percentage}\\
\hline
 1 & 1\%\\
\hline
 2 & 41\%\\
\hline
 3 & 45\%\\
\hline
 4 & 7\%\\
\hline
 unproven & 6\%\\
\hline
\end{tabular}
\label{tab:4}
\end{table}
More interesting is the plot of curves of rank 3 (Figure \ref{fig:2})
which shows that progression for curves of rank 3 is almost linear. It suggests that we can use the general algorithm from the introduction to state another version of the main theorem \ref{main_theorem} and find many more infinite families of elliptic curves over $\mathbb{Q}$. 
\subsection{Explicit version of the main theorem}
The statement of the theorem requires removing a finite subset of 'bad rational points' due to Falting's Theorem(Mordell's Conjecture). The upper bound of heights of this points is hard to obtain. We shall give an explicit and effective version of the main result of this paper.  However for rational points on the curve (\ref{explicit_eq}) with low height we can compute the explicit table of corresponding elliptic curves $E_u$ over $\mathbb{Q}$ of rank at least 4.
The curve:
\begin{equation}
C_1:y_0^2 = x_0^3 - 92835 x_0 + 1389150
\end{equation}
is mapped to the curve:
\begin{equation}
C_2:2569 + 18 u - 9 u^2 - 18 u^3 + 9 u^4=v^2
\end{equation}
via the map:
$$\phi:C_1\longrightarrow C_2$$
where $\phi(x_0,y_0)=(u,v)$ is given by the formulas:
\begin{equation*}
u=\frac{565605+x_0 (-948+7 x_0)+266 y_0}{(-1551+x_0) (45+x_0)}
\end{equation*}
\begin{align*}
v=\frac{1}{(-1551+x_0)^2 (45+x_0)^2}\left(133 (92385+(-30+x_0) x_0) (-115425+x_0 (1536+x_0))\right.+\\
\left. 234 (-3922935+x_0 (9010+41 x_0)) y_0\right)
\end{align*}
and defined at the points $(-45, 2340)$,$(-45,-2340)$, $(1551,59904)$, $(1551,-5990)$, $\infty_{C_1}$:
\begin{align*}
\phi(-45, 2340)&=\infty_{C_2}\\
\phi(-45,-2340)&=(-\frac{10898}{5187},-\frac{477412081}{8968323})\\
\phi(1551,59904)&=\infty_{C_2}\\
\phi(1551,-59904)&=(\frac{16085}{5187},\frac{477412081}{8968323})\\
\phi(\infty_{C_1})&=(7,133)\\
\end{align*}
for $\infty_{C_1}$ - the point at infinity $C_1$ and analogously for $C_2$. The map is regular at every point of $C_1$ so it is a morphism of curves. For the converse we have the mapping:
$$\psi:C_1\longrightarrow C_2$$
where $\psi(u,v)=(x_0,y_0)$ with:
$$x_0=\frac{5117-948 u+753 u^2+266 v}{(-7+u)^2}$$
$$y_0=\frac{266 (5201+9 u (-4+u (-22+13 u)))+2 (1799+4797 u) v}{(-7+u)^3}$$
which is not regular at the point $\infty_{C_2}$ and is defined at the points $(7,133)$ and $(7,-133)$:
\begin{align*}
\psi(7,133)&=\infty_{C_1}\\
\psi(7,-133)&=(-\frac{3628425}{17689},\frac{8081948160}{2352637})
\end{align*}
If we define sets:
$$A=C_1\backslash\{(-45, 2340),(1551,59904)\}$$
$$B=C_2\backslash\{\infty_{C_2}\}$$
Then we have:
$$\phi\circ\psi=id_A$$
$$\psi\circ\phi=id_B$$
We now give the explicit table of curves of rank at least 4 as stated in the main theorem. If we assume the Parity Conjecture we can show that some of them have actually the rank at least 5. Let $E_u:y^2+(u^2-u-3)xy=x^3+(u^2-u-3)x^2-x+1$ and $P_1=(-309, 756)$, $P_2=(-45,2340)$, $T=(15,0)$ - the points spanning the group $C_1(\mathbb{Q})$. From the computations above we can associate uniquely a pair $(u,v)$ on $C_2$ corresponding to the point $\alpha T+\beta_1 P_1+\beta_2 P_2$. We abbreviate this as $(u,v)\leftrightarrow(\alpha,\beta_1,\beta_2)$. We define the following functions:
\begin{itemize}
\item $R(u)$ is the regulator of points $(0,1),(1,1),(u,u+1),\left(\frac{1}{9},\frac{1}{54} (9 + 3 u - 3 u^2 + v)\right)$
\item $N(u)$ is the conductor of the curve $E_u$
\item $j(u)$ is the j-invariant of $E_u$
\item $W(u)$ is equal to the global root number $W(E_u/\mathbb{Q})$.
\end{itemize}
All the computations were performed for the minimal model of each curve.
\begin{table}[h]
\caption{Curves of rank 4 and 5}
\centering
\makebox[0cm][c]{
\begin{tabular}{|c|c|c|c|c|c|}
\hline
$(\alpha,\beta_1,\beta_2)$ & $\textbf{R(u)}$ & $\textbf{N(u)}$ &$\textbf{j(u)}$ & $\textbf{W(u)}$ & \textbf{Rank}\\
\hline
( 0 , -2 , -2 ) & 253637.08 & $7.42 \times 10^{117}$ & -4382.17 & -1 & $\geq 5$\\
\hline
( 0 , -2 , -1 ) & 53400.57 & $4.79 \times 10^{79}$ & $-1.39 \times 10^{6}$ & 1 &$\geq 4$\\
\hline
( 0 , -2 , 0 ) & 16681.20 & $5.69 \times 10^{59}$ & $-4.14\times 10^{11}$ & -1 & $\geq 5$\\
\hline
( 0 , -2 , 1 ) & 23528.39 & $1.89 \times 10^{64}$ & $-1.11 \times 10^{16}$ & 1 & $\geq 4$\\ 
\hline
( 0 , -1 , -2 ) & 117347.77 & $1.22 \times 10^{95}$ & $-4.66 \times 10^{19}$ & 1  & $\geq 4$\\
\hline
( 0 , -1 , -1 ) & 6398.35 & $2.46 \times 10^{46}$ & $-7.42 \times 10^{8}$ & 1 & $\geq 4$\\
\hline
( 0 , -1 , 0 ) & 28.40 & $4.13 \times 10^{12}$ & -1255.79 & 1 & $\geq 4$\\
\hline
( 0 , -1 , 1 ) & 0.0 & $6.54 \times 10^{11}$ & -1264.95 & 1 & -\\
\hline
( 0 , 0 , -2 ) & 138113.04 & $6.46 \times 10^{98}$ & -912.11 & -1 & $\geq 5$\\
\hline
( 0 , 0 , -1 ) & 4697.68 & $1.21 \times 10^{43}$ & -20742.18 & 1 & $\geq 4$\\
\hline
( 0 , 0 , 0 ) & 8.61 & 57482738.0 & $-4.72 \times 10^{9} $ & 1 & $\geq 4$\\
\hline
( 0 , 1 , -2 ) & 608830.99 & $3.64 \times 10^{145}$ & $-4.45 \times 10^{19}$ & 1 & $\geq 4$\\
\hline
( 0 , 1 , -1 ) & 56796.71 & $1.80 \times 10^{81}$ & $-3.75 \times 10^{10}$ & -1 & $\geq 5$\\
\hline
( 0 , 1 , 0 ) & 1301.45 & $8.98 \times 10^{31}$ & -200862.89 & -1  & $\geq 5$\\
\hline
( 0 , 1 , 1 ) & 0.0 & $6.54 \times 10^{11}$ & -1264.95 & 1 & -\\
\hline
\end{tabular}
}
\label{tab:5}
\end{table}
For the last tuple the regulator is equal to 0 because the tuple corresponds to $u=\frac{1}{9}$ for which the fourth point from the statement of Theorem \ref{main_theorem} coincides with the third point. For the tuple $(0,-1,1)$ (when $u=8/9$) the fourth point is linearly dependent on the other three points. Moreover the curves corresponding to these tuples are isomorphic over $\mathbb{Q}$.
\begin{remark}
We can find in the family $E_u$ curves of unconditional rank at least five. The curve $E_{16}:y^2+239 x y=x^3+239 x^2-x+1$ is a curve of unconditional rank five. The set of generators of the non-torsion part is given by:
$$(0,1),(1,1),(16,17),(-14/25,16661/125),(52/81,469/729)$$
We can show that for $c=-14/25$ the associated auxiliary elliptic curve from Theorem \ref{main_theorem}:
$$4 - 4 c - 3 c^2 + 4 c^3 + 2 c^2 u - c^2 u^2 - 2 c^2 u^3 + c^2 u^4=t^2$$
has rank 4 over $\mathbb{Q}$.
This might lead to alternative formulation of Theorem \ref{main_theorem}.
\end{remark}
\section*{Acknowledgments}
The author would like to thank Wojciech Gajda for suggesting the problem. He thanks Sebastian Petersen for helpful comments concerning root numbers and the Parity Conjecture. The author is grateful to Adam Lipowski for the computational resources supporting SAGE 3.4.
\bibliographystyle{amsplain}

\end{document}